\documentclass[a4paper,12pt]{article} %izvzela 'twoside' fleqn,

\usepackage{color}

\usepackage{amsfonts, amsmath, amsthm, amssymb}
\usepackage[T1]{fontenc}
\usepackage[cp1250]{inputenc}
\usepackage{xcolor}
\usepackage{graphicx}
\usepackage{amssymb}
\usepackage{amsmath}
\usepackage{mathptmx}
\usepackage{helvet}
\usepackage{courier}
\usepackage{txfonts}
\usepackage{tikz} 
\usetikzlibrary{arrows}
\usepackage{type1cm}

\usepackage{verbatim}

\usepackage{graphicx}
\usepackage{epsfig,amscd,amssymb,amsxtra,amsmath,amsthm}
\usepackage{type1cm}
\usepackage[T1]{fontenc}
\usepackage{graphics}
\usepackage[mathscr]{eucal}
\usepackage[all]{xy}
\usepackage{amsmath,amscd}

\newtheorem{theorem}{Theorem}[section]

\newtheorem{definition}[theorem]{Definition}
\newtheorem{lemma}[theorem]{Lemma}

\newtheorem{example}[theorem]{Example}

\newtheorem{corollary}[theorem]{Corollary}
\newtheorem{problem}[theorem]{Problem}
\newtheorem{observation}[theorem]{Observation}

\newcommand{\Bd}  {\mathop{\rm Bd}\nolimits}

\begin{document}

\def\joinrel{\mkern-3mu}
\newcommand{\varproj}{\displaystyle \lim_{\multimapinv\joinrel-\joinrel-}}

\title{Towards the complete classification of fans}
\author{I. Bani\v c, G. Erceg, A. Illanes, I. Jeli\' c, J. Kennedy, V. Nall}
\date{}

\maketitle

\begin{abstract}
A fan is an arcwise-connected continuum, which is hereditarily unicoherent and has exactly one ramification point. Many of the known examples of fans were constructed as 1-dimensional continua that are unions of arcs which intersect in exactly one point. Borsuk proved in 1954 that each fan is a $1-$dimensional continuum which is the union of arcs intersecting in exactly one point. But it is not yet known if this property is equivalent to being a fan. In this paper, we show that under an additional assumption, every such union of arcs is a fan.
\end{abstract}
\-

%%%%%%%%%%%%%%%%%%%%%%%%%%%%%%%%%%%%%%%%%%%%%%%%%%%%%%%%%%%%%%%%%%%%%%%%%%%%%%%%%
%%% I N T R O D U C T I O N S
\section{Introduction}
Fans and their properties (smoothness, property of Kelley \cite{charatonik,charatonik1}, and others) have been studied  intensively. They are also very important when studying hyperspaces (for examples see \cite{nadler,23,macias}).  Many interesting examples of fans have been constructed and studied. For example, the Borsuk fan, which can not be embedded into the plane \cite{borsuk}; the Lelek fan, a fan with $1-$dimensional set of end-points \cite{charatonik3,lelek}; the Cantor fan \cite{charatonik2}; and others \cite{jjc,jjc1,jjc2,jjc3,charatonik,charatonik1}. 
Most of the examples of fans were constructed as unions of arcs which intersect in exactly one point. The fact that a fan $X$ with the set of end-points $E(X)$ can be written as the union of arcs joining end-points to the top of $X$, where the union is taken over $E(X)$ was already mentioned in \cite{jjc} as something that was easy to see. In fact, Theorem \ref{c1} is a well-known result. It was proved by Borsuk in \cite{borsuk1}.

\begin{theorem} \label{c1} Let $X$ be a fan with the top $v$. Then there is a family of arcs $\mathcal L$ in $X$ such that
\begin{enumerate}
\item\label{unk} $\displaystyle X=\bigcup_{L\in \mathcal L}L$;
\item\label{dunk} for all $L_1,L_2\in \mathcal L$,
$$
L_1\neq L_2 ~~~ \Longrightarrow ~~~  L_{1}\cap L_{2}=\{v\}.
$$
\end{enumerate}
\end{theorem}
Example \ref{tripikaena} shows that for a continuum $X$, properties \ref{unk} and \ref{dunk} from Theorem \ref{c1} do not imply that $X$ is a fan. But it is not yet known if $X$ is a $1$-dimensional continuum, whether these two properties are equivalent to being a fan. 

We proceed as follows. In Section \ref{s2}, we give basic definitions and results that are used later. In Section \ref{s3}, we discuss the converse of Theorem \ref{c1}, and, in Section \ref{s4}, we provide a large family of $1$-dimensional continua which are the unions of arcs intersecting in exactly one point 
that are fans. 

%We still don't know if each $1$-dimensional continuum which is the union of arcs intersecting in exactly one point is a Carolyn phan; see Problem \ref{prlekija12}.
	\section{Definitions and notation}\label{s2}
	In this section, we give basic definitions that are used later in the paper.
	\begin{definition}
 \emph{A continuum} is a non-empty compact connected metric space.  \emph{A subcontinuum} is a subspace of a continuum, which is itself a continuum.
 \end{definition}
\begin{definition}\label{ccc}
Let $X$ be a continuum. 
\begin{enumerate}
\item The continuum $X$ is \emph{unicoherent} if for any subcontinua $A$ and $B$ of $X$ such that $X=A\cup B$,  the compactum $A\cap B$ is connected. 
\item The continuum $X$ is \emph{hereditarily unicoherent } provided that each of its subcontinua is unicoherent.
\item The continuum $X$ is a \emph{dendroid} if it is an arcwise connected, hereditarily unicoherent continuum.
\item If $X$ is homeomorphic to $[0,1]$, then $X$ is \emph{an arc}.   
\item Let $X$ be an arc. A point $x\in X$ is called \emph{an end-point of the arc  $X$} if  there is a homeomorphism $\varphi:[0,1]\rightarrow X$ such that $\varphi(0)=x$.
\item Let $X$ be a dendroid.  A point $x\in X$ is called an \emph{end-point of the dendroid $X$} if for  every arc $A$ in $X$ that contains $x$, $x$ is an end-point of $A$.  The set of all end-points of $X$ is denoted by $E(X)$. 
\item The continuum $X$ is \emph{a simple triod} if it is homeomorphic to $([-1,1]\times \{0\})\cup (\{0\}\times [0,1])$.
\item Let $X$ be a simple triod. A point $x\in X$ is called \emph{the top-point} or, briefly, the \emph{top of the simple triod $X$} if  there is a homeomorphism $\varphi:([-1,1]\times \{0\})\cup (\{0\}\times [0,1])\rightarrow X$ such that $\varphi(0,0)=x$.
\item Let $X$ be a dendroid.  A point $x\in X$ is called \emph{a ramification-point of the dendroid $X$}, if there is a simple triod $T$ in $X$ with top $x$.  The set of all ramification-points of $X$ is denoted by $R(X)$. 
\item The continuum $X$ is \emph{a fan} if it is a dendroid with at most one ramification point $v$, which is called the top of the fan $X$ (if it exists).
\item Let $X$ be a fan.   For all points $x$ and $y$ in $X$, we define  \emph{$[x,y]$} to be the arc in $X$ with end-points $x$ and $y$, if $x\neq y$. If $x=y$, then we define $[x,y]=\{x\}$.
\item Let $X$ be a fan with top $v$. We say that that the fan $X$ is \emph{smooth} if for any $x\in X$ and for any sequence $(x_n)$ of points in $X$,
$$
\lim_{n\to \infty}x_n=x \Longrightarrow \lim_{n\to \infty}[v,x_n]=[v,x].
$$ 
A fan is non-smooth if it is not smooth. 
\item Let $X$ be a fan with the top $v$. Then we use $\mathcal L(X)$ to denote the family of arcs in $X$ such that
\begin{enumerate}
\item $\displaystyle X=\bigcup_{L\in \mathcal L(X)}L$;
\item for all $L_1,L_2\in \mathcal L(X)$,
$$
L_1\neq L_2 ~~~ \Longrightarrow ~~~  L_{1}\cap L_{2}=\{v\}.
$$
\end{enumerate}
The arcs in $\mathcal L(X)$ are called legs of the fan $X$.
\item Let $X$ be a fan.  We say that $X$ is \emph{a Cantor fan} if $X$ is homeomorphic to the continuum $\bigcup_{c\in C}A_c$, where $C\subseteq [0,1]$ is the standard Cantor set and for each $c\in C$, $A_c$ is the  {convex} segment in the plane from $(\frac{1}{2},0)$ to $(c,1)$.
%\item Let $X$ be a fan.  We say that $X$ is \emph{a Lelek fan} if it is smooth and $\Cl(E(X))=X$. See Figure \ref{figure2}.
%\begin{figure}[h!]
%	\centering
%		\includegraphics[width=20em]{Lelek_FIG.pdf}
%	\caption{A Lelek fan}
%	\label{figure2}
%\end{figure}  
\end{enumerate}
\end{definition}
\begin{observation} \label{tatoo}
It is a well-known fact that the Cantor fan is universal for smooth fans, {i.e., every smooth fan embeds into it} (for details see \cite[Theorem 9, p. 27]{jjc},  \cite[Corollary 4]{koch},  and \cite{eberhart}). %Since the Lelek fan contains a Cantor fan (for details see \cite{short}), it follows  also that the Lelek fan is a universal continuum for smooth fans.  
\end{observation}
\begin{definition}
Let $(X,d)$ be a compact metric space. Then we define \emph{$2^X$} by 
$$
2^{X}=\{A\subseteq X \ | \ A \textup{ is a non-empty closed subset of } X\}.
$$
Let $\varepsilon >0$ and let $A\in 2^X$. Then we define  \emph{$N_d(\varepsilon,A)$} by 
$$
N_d(\varepsilon,A)=\bigcup_{a\in A}B(a,\varepsilon).
$$
The function \emph{$H_d:2^X\times 2^X\rightarrow \mathbb R$}, defined by
$$
H_d(A,B)=\inf\{\varepsilon>0 \ | \ A\subseteq N_d(\varepsilon,B), B\subseteq N_d(\varepsilon,A)\},
$$
for all  $A,B\in 2^X$, is called \emph{the Hausdorff metric} on $2^X$.  
\end{definition}
Let $(X,d)$ be a compact metric space. The Hausdorff metric on $2^X$ is in fact a metric on $2^X$ and the metric space $(2^X,H_d)$ is called \emph{the hyperspace of the space $(X,d)$}. Also, let $A$ be a non-empty closed subset of $X$,  and let $(A_n)$ be a sequence of non-empty closed subsets of $X$. When we write $\displaystyle A=\lim_{n\to \infty}A_n$, we mean $\displaystyle A=\lim_{n\to \infty}A_n$ in $(2^X,H_d)$. 

\begin{definition}
	Let $(X,d)$ be a compact metric space, let $A$ be a non-empty closed subset of $X$,  and let $(A_n)$ be a sequence of non-empty closed subsets of $X$. We use $\limsup A_n$ to denote the set of points $x\in X$ such that for each open set $U$ in $X$ that contains $x$, $U\cap A_n\neq \emptyset$ for infinitely many positive integers $n$.
\end{definition}

We use $\dim (X)$ to refer to the small inductive dimension \cite[Definition 1.1.1., page 3]{nagata}. It is also called the (topological) dimension of topological spaces.  We use the following well-known result about $1$-dimensional continua, which follows directly from \cite[Definition 1.1.1., page 3]{nagata}.
\begin{theorem}\label{enadim}
	Let $X$ be a $1$-dimensional continuum. Then for each $x\in X$ and for each open neighborhood $U$ of $x$ in $X$, there is an open  neighborhood $V$ of $x$ in $X$ such that
	\begin{enumerate}
		\item $x\in V\subseteq U$,
		\item $\dim (\Bd(V))=0$. 
	\end{enumerate}
	(We use $\Bd$ to denote the boundary of a set in a topological space.)
\end{theorem}
\noindent We also use the following well known result, which follows from \cite[Theorem 1.4.5., page 33]{nagata}.
\begin{theorem}\label{totally}
	Let $X$ be a compact metric space. The following statements are equivalent.
	\begin{enumerate}
		\item $\dim (X)=0$. 
		\item $X$ is totally disconnected.
	\end{enumerate}
\end{theorem}
\section{The converse of Theorem \ref{c1}}\label{s3}
In this section, we discuss the converse of Theorem \ref{c1} and prove the converse of Theorem \ref{c1} under the assumption that the continuum is hereditarily unicoherent.  

\begin{figure}[h!]
 \begin{center}
    \includegraphics[width=4in]{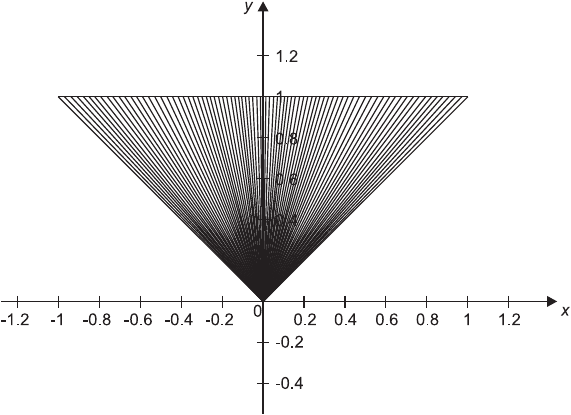}
    \caption{$X=\bigcup_{x\in A}L_{x}$}
    \label{fan}
  \end{center}
\end{figure}

The following example shows that, if for a continuum $X$ there is a family of arcs $\mathcal L$ in $X$ such that
\begin{enumerate}
\item $\displaystyle X=\bigcup_{L\in \mathcal L}L$;
\item for all $L_1,L_2\in \mathcal L$,
$$
L_1\neq L_2 ~~~ \Longrightarrow ~~~  L_{1}\cap L_{2}=\{v\}.
$$
\end{enumerate}
then $X$ may not be a fan.
\begin{example}\label{tripikaena}
Let $A=[-1,1]\times \{1\}\subseteq \mathbb{R}^2$, and $X=\bigcup_{x\in A}L_{x}$, where for each $x\in A$, $L_{x}$ is the line segment from $(0,0)$ to $x$ in $\mathbb{R}^2$, see Figure \ref{fan}. Obviously $X$ is not a fan, as it is not hereditarily unicoherent (note that it contains a simple closed curve).
\end{example}
%%%%%%%%%%%%%%%%%%%%%%%%
%% SLIKA 
%%%%%%%%%%%%%%%%%%%%%%%%%%
%%%%%%%%%%%%%%%%%%%%%POGLAVJE
%%%%%%%%%%%%%%%%%%%%%%%%%%%%%%%%%%%%%%%%%%%%%%%%%%%%%%%%%%%%%%%%%%%%%%%%%%%%%%%%%%%%%%%%%%%%%%%%%%%%%%%%%%%%%%%%%%%%
All fans are dendroids and all dendroids are $1$-dimensional continua (\cite[(48), p.\ 239]{jjc1}). Therefore, to prove the converse of Theorem \ref{c1} it suffices to consider only 1-dimensional continua. 
Next, we prove Lemma \ref{usledir}.
\begin{lemma}\label{usledir}
Let $X$ be a  hereditarily unicoherent continuum and $v\in X$. Suppose there is a family of arcs $\mathcal L$ in $X$ such that $|\mathcal L|\geq 3$ and
\begin{enumerate}
\item $\displaystyle X=\bigcup_{L\in \mathcal L}L$;
\item for all $L_1,L_2\in \mathcal L$,
$$
L_1\neq L_2 ~~~ \Longrightarrow ~~~  L_{1}\cap L_{2}=\{v\}.
$$
\end{enumerate}
Then $v$ is the only ramification point in $X$.
\end{lemma}
\begin{proof}
Let $a\in X\setminus \{v\}$ be another ramification point in $X$. Then there is a triod $T$ in $X$ with top $a$. Let $A_1,A_2,A_3$ be arcs in $T$, such that $T=A_1\cup A_2\cup A_3$ and $A_1\cap A_2=A_2\cap A_3=A_1\cap A_3=\{a\}$. Next, let $L_0\in \mathcal L$ be such that $a\in L_{0}$. Since $L_{0}$ is an arc, one of the arcs $A_1,A_2,A_3$ intersects with $L_{0}$ only at the point $a$. Without loss of generality, suppose $A_1\cap L_{0}=\{a\}$. If there is $L_1\in \mathcal L$, such that $A_1\subseteq L_{1}$, then $L_{0}\cap L_{1}=\{v,a\}$, which contradicts the assumptions of this lemma. Therefore, there are $x,y\in A_1$, $x\neq y$, and $L_1, L_2\in \mathcal L$, $L_1\neq L_2$, such that $x\in A_1\cap L_{1}$ and $y\in A_1\cap L_{2}$. In this case $A_1\cup L_{1}\cup L_{2}$ is a subcontinuum of $X$, which is not unicoherent. Therefore, $v$ is the only ramification point in $X$.    
\end{proof}
Finally, the Theorem \ref{main2} follows directly from Lemma \ref{usledir}.
\begin{theorem}\label{main2}
Let $X$ be a $1$-dimensional hereditarily unicoherent continuum and $v\in X$. Suppose that there is a family of arcs $\mathcal L$ in $X$ such that $|\mathcal L|\geq 3$ and
\begin{enumerate}
\item $\displaystyle X=\bigcup_{L\in \mathcal L}L$;
\item for all $L_1,L_2\in \mathcal L$,
$$
L_1\neq L_2 ~~~ \Longrightarrow ~~~  L_{1}\cap L_{2}=\{v\}.
$$
\end{enumerate}
Then $X$ is a fan.
\end{theorem}

We conclude this section by stating the following open problem.
\begin{problem}\label{main3}
Let $X$ be a $1$-dimensional continuum and $v\in X$. Are the following statements equivalent?
\begin{enumerate}
\item $X$ is a fan with top $v$.
\item There is a family  $\mathcal L$ of arcs in $X$ such that $|\mathcal L|\geq 3$ and 
\begin{enumerate}
\item $\displaystyle X=\bigcup_{L\in \mathcal L}L$;
\item for all $L_1,L_2\in \mathcal L$,
$$
L_1\neq L_2 ~~~ \Longrightarrow ~~~  L_{1}\cap L_{2}=\{v\}.
$$
\end{enumerate}
\end{enumerate}
\end{problem}

\section{Pans and phans versus fans}\label{s4}
In this section, we provide a large family of $1$-dimensional continua which are the unions of arcs intersecting in exactly one point (we call it a family of Matea phans) that are fans.

\begin{definition}
	Let $X$ be a continuum, let $v\in X$ and let $\mathcal L$ be a collection of arcs in $X$. We say that $(v,\mathcal L)$ is \emph{a ramification pair in $X$} if 
	\begin{enumerate}
	\item $\displaystyle \bigcup_{L\in \mathcal L}L=X$, and
		\item for all $L_1,L_2\in \mathcal L$,
		$$
		L_1\neq L_2 ~~~ \Longrightarrow ~~~ L_{1}\cap L_{2}=\{v\}.
		$$
	\end{enumerate}
\end{definition}

\begin{definition}
	Let $X$ be a continuum and let $(v,\mathcal L)$ be a ramification pair in $X$. Then we say that $X$ is \emph{a pan with respect to $(v,\mathcal L)$}. We say that $X$ is \emph{ a pan} if  there a ramification pair $(v,\mathcal L)$ in $X$ such that $X$ is a pan with respect to $(v,\mathcal L)$. 
\end{definition}

\begin{definition}
	Let $X$ be a pan, let $v\in X$ and let $\mathcal L$ be a collection of arcs in $X$. We say that $X$ is a \emph{Carolyn pan with respect to $(v,\mathcal L)$} if $X$ is a pan with respect to $(v,\mathcal L)$ and for each continuum $C$ in $X$, 
	$$
	v\not \in C ~~~  \Longrightarrow  ~~~ \textup{there is } L\in \mathcal L \textup{ such that } C\subseteq L.
	$$
	We say that $X$ is \emph{a Carolyn pan} if there are a point $v\in X$ and a family $\mathcal L$ of arcs in $X$ such that $X$ is a Carolyn pan with respect to $(v,\mathcal L)$. 
%	We say that $X$ is a Carolyn phan, if $X$ is a phan as well as a Carolyn pan. If $X$ is a Carolyn phan as well as a pan with respect to $(v,\mathcal L)$, then we say that $X$ is a Carolyn phan with respect to $(v,\mathcal L)$
\end{definition}

\begin{observation}
	Note that each fan is a Carolyn pan.
\end{observation}
	\begin{observation}
Example \ref{tripikaena} is a pan which is not a Carolyn pan with respect to $\big((0,0), \{L_x \ | \ x\in A\}\big)$. Therefore, the topological product $\prod_{i=1}^2[0,1]$ is a pan. One can use a similar argument to prove that for each $n\in \{\infty, 1,2,3,\ldots\}$, the topological product $\prod_{i=1}^n[0,1]$ is a pan. It follows that for each $n\in \{\infty, 1,2,3,\ldots\}$, there is a pan $X$ such that $\dim(X)=n$. Moreover, note that for any compactum $X$, the cone over $X$ is a Carolyn pan. 
\end{observation}

\begin{definition}
	Let $X$ be a pan. If $\dim (X)=1$, we say that $X$ is \emph{a phan}. 
	\end{definition}
	\begin{definition}
	Let $X$ be a phan, let $v\in X$ and let $\mathcal L$ be a collection of arcs in $X$.
	We say that $X$ is \emph{a phan  with respect to $(v,\mathcal L)$} if $X$ is a pan with respect to $(v,\mathcal L)$.
		\end{definition}

\begin{definition}
We use $S$ to denote the unit circle 
$$
S=\{(x,y)\in \mathbb R^2 \mid x^2+z^2=1\}
$$
 in the plane $\mathbb{R}^{2}$. We also use 
\begin{center}
$S^{+}=\{(x,y)\in S\mid y\geq 0\}$, and

$S^{-}=\{(x,y)\in S\mid y\leq 0\}$.
\end{center}
We also use $e$ to denote the exponential mapping $e:\mathbb{R}\rightarrow S$ given by
\begin{center}
$e(t)=(\operatorname{cos}(2\pi t),\operatorname{sin}(2\pi t))$
\end{center}
for each $t\in \mathbb R$. 
\end{definition}
\begin{definition}
	Let $X$ be a pan with respect to $(v,\mathcal L)$ and let $x,y\in X$. If $x\neq y$ and there is $L\in \mathcal L$ such that $x,y\in L$, then we use $[x,y]_{(v,\mathcal L)}$ to denote the arc in $L$ with end-points $x$ and $y$. If $x=y$, then we define $[x,y]_{(v,\mathcal L)}=\{x\}$.
\end{definition}
\begin{definition}
	Let $X$ be a phan with respect to $(v,\mathcal{L})$. For each $p\in X\setminus\{v\}$, we use  $L_{p}$ to denote the unique element of $\mathcal{L}$ containing $p$. 
	\end{definition}
	\begin{definition}
	Given a subcontinuum $B$ of $X$ we say that $B$ is $(v,\mathcal{L})$-{\it star} if for each $b\in B$, $[v,b]_{(v,\mathcal{L})}\subset B$.
\end{definition}
\begin{definition}
A continuum $X$ is said to be \emph{contractible} if there exists a point 
$x_0 \in X$ and a continuous map 
\[
H : X \times [0,1] \to X
\]
such that for every $x \in X$:
\[
H(x,0) = x \quad \text{and} \quad H(x,1) = x_0.
\]
That is, the identity map on $X$ is homotopic to a constant map.
\end{definition}
\begin{definition}
We say that $X$ is a {\it Matea} phan with respect to $(v,\mathcal{L})$ if for each sequence $(L_{n})$ in $\mathcal{L}$, convergent to a subcontinuum $B$ of $X$, we have that $B$ is contractible and $(v,\mathcal{L})$-star. 
\end{definition}

\begin{theorem} \label{enica}
Every Matea phan is a Carolyn pan.
\end{theorem}

\begin{proof}
Let $X$ be a Mathea phan with respect to $(v,\mathcal{L})$. Suppose that $X$ is not Carolyn. Then there exists a subcontinuum $C$ of $X$ such that $v\notin C$ and $C$ is not contained in any arc in $\mathcal{L}$.
Thus there exist distinct arcs $J$ and $K$ in $\mathcal{L}$ and points $p_{J}\in J\cap C$ and $p_{K}\in K\cap C$.

By Urysohn Lemma there exist mappings $f_{J}:J\cup C\rightarrow S^{+}$, $f_{K}:K\cup C\rightarrow S^{-}$ such that $f_{J}(v)=(1,0)=f_{K}(v)$ and $f_{J}(C)=\{(-1,0)\}=f_{K}(C)$. Since $(J\cup C)\cap (K\cup C)=\{v\}\cup C$, there exists a mapping $f^{*}:J\cup C\cup K\rightarrow S$ that extends $f_{J}$ and $f_{K}$.

Since $\operatorname{dim}[X]=1$, by \cite[Theorem VI 4, p. 83]{hurewicz}, there exists a mapping $f:X\rightarrow S$ that extends $f^{*}$. Given $p\in C$, since $p\neq v$, $L_{p}$ is defined. By \cite[Theorem 12.9, p. 155]{Illanes2025}, there exists a lifting $h_{p}:L_{p}\rightarrow \mathbb{R}$ such that $f\vert_{L_{p}}=e\circ h_{p}$. Observe that $f(v)=(1,0)$ and $f(C)=\{(-1,0)\}$. Then we can assume that $h_{p}(v)=0$, and we have that $h_{p}(p)$ is of the form $h_{p}(p)=j_{p}+\frac{1}{2}$, where $j_{p}$ is an integer, observe that $j_{p}$ is unique.

Let
\begin{center}
$C^{+}=\{p\in C: j_{p}\geq 0\}$ and $C^{-}=\{p\in C: j_{p}<0\}$.
\end{center}

Observe that $e\vert_{[0,\frac{1}{2}]}:[0,\frac{1}{2}]\rightarrow S^{+}$ is a homeomorphism satisfying $e\vert_{[0,\frac{1}{2}]}(0)=(1,0)$. Since $f(v)=(1,0)$, we obtain that the unique lifting $h:J\cup C\rightarrow \mathbb{R}$ of $f\vert_{J\cup C}=f_{J}$ satisfying $h(v)=0$ is the mapping $(e\vert_{[0,\frac{1}{2}]})^{-1}\circ f_{J}$. Thus $h_{p_{J}}=((e\vert_{[0,\frac{1}{2}]})^{-1}\circ
f_{J})\vert_{J}$. Hence $h_{p_{J}}(p_{J})=\frac{1}{2}$. We have shown that $j_{p_{J}}=0$ and then $C^{+}\neq\emptyset$. Similarly, using the mapping $e\vert_{[-\frac{1}{2},0]}$, we obtain that $j_{p_{K}}=-1$, so $C^{-}\neq\emptyset$.

Since $C=C^{+}\cup C^{-}$ and $C$ is connected, we have that $C^{+}$ and $C^{-}$ are not separated. Then we can suppose that there exists a point $p\in C^{+}\cap \operatorname{cl}_{X}(C^{-})$ the case $C^{-}\cap \operatorname{cl}_{X}(C^{+}) \neq \emptyset$ is similar. 
Let $(p_{n})$ be a sequence in $C^{-}$ converging to $p$. Taking a subsequence, if necessary, we may assume that the sequence $(L_{p_n})$ converges to a subcontinuum $B$ of $X$. Then $p\in B$. Since $X$ is a Matea phan, $B$ is contractible and $(v,\mathcal{L})$-star.

Since $B$ is contractible, by \cite[Theorem 12.7, p. 151]{Illanes2025}, the mapping $f\vert_{B}:B\rightarrow S$ can be lifted. By \cite[Theorem 12.6, p. 150]{Illanes2025}, there exists an open subset $U$ of $X$ and a lifting $h_{U}:U\rightarrow \mathbb{R}$ of $f\vert_{U}$ such that $B\subset U$, we may assume that $h_{U}(v)=0$. 

By the continuity of $h_{U}$, there exists $m\in\mathbb{N}$ such that $\vert h_{U}(p)-h_{U}(p_{m})\vert<\frac{1}{4}$ and since $(L_{p_m})$ converges to $B$, we may suppose that $L_{p_{m}}\subset U$. Then $h_{U}\vert_{L_{p_{m}}}$ and $h_{p_{m}}$ are liftings of $f_{L_{p_{m}}}$ that coincide in $v$, so $h_{U}\vert_{L_{p_{m}}}=h_{p_{m}}$ and $h_{U}(p_{m})=h_{p_{m}}(p_{m})=j_{p_{m}}+\frac{1}{2}$. Since $B$ is $(v,\mathcal{L})$-star, we have that $[v,p_{}]_{(v,\mathcal{L})}\subset B\subset U$. This implies that $h_{U}\vert_{[v,p_{}]_{(v,\mathcal{L})}}$ and $h_{p}\vert_{[v,p_{}]_{(v,\mathcal{L})}}$ are liftings of $f\vert_{[v,p_{}]_{(v,\mathcal{L})}}$ that coincide in $v$. Thus $h_{U}(p)=h_{p}(p)=j_{p}+\frac{1}{2}$. Therefore $\vert j_{p}+\frac{1}{2}-(j_{p_{m}}+\frac{1}{2})\vert<\frac{1}{2}$, so $j_{p}=j_{p_{m}}$. This contradicts the fact that $j_{p}\geq 0$ and $j_{p_{m}}<0$. Therefore $X$ is Carolyn.  \end{proof}

\begin{theorem}\label{dvojica}
	Let $X$ be a phan. If $X$ is a Carolyn pan, then $X$ is a fan.
	\end{theorem}
\begin{proof}
	Let $X$ be a Carolyn pan with respect to $(v,\mathcal L)$. First, we prove that $X$ is a hereditarily unicoherent continuum.  Suppose that $X$ is not hereditarily unicoherent. Let $A$, $B$ and $C$ be continua in $X$ such that $C=A\cup B$ and $A\cap B$ is not connected. Also, let $S,T$ be a separation of $A\cap B$. Suppose that $v\not\in C$. It follows that for some $L\in \mathcal L$, $C\subseteq L$ (since $X$ is a Carolyn pan). Therefore, $C$ is an arc, which is a contradiction since $C$ is not unicoherent. Thus, $v\in C$. Suppose that $v\in A\setminus B$ (or $v\in B\setminus A$).  Let $x\in S$ and $y\in T$. Since $v\not\in B$ it follows that there is $L\in \mathcal L$ such that $B\subseteq L$ (such an $L$ does exist since $X$ is a Carolyn pan). Choose and fix such an $L$. Let $x_0=x$ and $y_0=y$.  For each positive integer $n$, let $x_n,y_n\in \Bd(B(v,\frac{d(x,y)}{2^n}))$ be such that for some continua $X_n$ and $Y_n$ in $B\setminus B(v,\frac{d(x,y)}{2^n})$, $x_n,x_{n-1}\in X_n$ and $y_n,y_{n-1}\in Y_n$. Such points and continua do exist by the Boundary bumping theorem \cite[Theorem 5.4]{nad}. For each positive integer $n$, choose and fix such points $x_n$ and $y_n$, and such continua $X_n$ and $Y_n$. Note that for each positive integer $n$, $v\not\in X_n\cup Y_n$. Therefore, for each positive integer $n$, $X_n\cup Y_n\subseteq L$. Since 	$\displaystyle\lim_{n\to \infty}x_n=\lim_{n\to \infty}y_n=v$, it follows that there is a simple closed curve in $L$, which is a contradiction. Therefore, $v\in A\cap B$. Suppose that $v\in S$ and let $x\in T$. Let $x_0=y_0=x$.  For each positive integer $n$, let $x_n,y_n\in \Bd(B(v,\frac{d(x,v)}{2^n}))$ be such that for some continua $X_n\subseteq A\setminus B(v,\frac{d(x,y)}{2^n})$ and $Y_n\subseteq B\setminus B(v,\frac{d(x,y)}{2^n})$, $x_n,x_{n-1}\in X_n$ and $y_n,y_{n-1}\in Y_n$. Such points and continua do exist by the Boundary bumping theorem \cite[Theorem 5.4]{nad}. For each positive integer $n$, choose and fix such points $x_n$ and $y_n$, and such continua $X_n$ and $Y_n$. Note that for each positive integer $n$, $v\not\in X_n\cup Y_n$. Therefore, there is $L\in \mathcal L$ such that for each positive integer $n$, $X_n\cup Y_n\subseteq L$. Choose and fix such an $L$. Since 	$\displaystyle\lim_{n\to \infty}x_n=\lim_{n\to \infty}y_n=v$, it follows that there is a simple closed curve in $L$, which is a contradiction. This proves that $X$ is a hereditarily unicoherent continuum. It follows from Lemma \ref{usledir} that $X$ is a fan.
	\end{proof}
	
\begin{corollary}
	Let $X$ be a continuum. If $X$ is a Matea phan, then $X$ is a fan.
	\end{corollary}
\begin{proof}
	 Suppose that $X$ is a Matea phan with respect to some $(v,\mathcal L)$. By Theorem \ref{enica}, $X$ is a Carolyn pan and by Theorem \ref{dvojica}, $X$ is a fan.
\end{proof}
We conclude the paper by stating the following open problem. 
	\begin{problem}\label{prlekija12}
	Is every phan a Carolyn pan?
\end{problem}

\section{Acknowledgement}
This work is supported in part by the Slovenian Research Agency (research projects J1-4632, BI-HR/23-24-011, BI-US/22-24-086 and BI-US/22-24-094, research program P1-0285) and University of Split, project code: IP-UNIST-44 (ITPEM).

%\section{Declarations}
%The following sections are not relevant to our manuscript. 
%\subsection{Competing interests}
%Not applicable.
%\subsection{Data Availability Statement}
%Not applicable.

\noindent I. Bani\v c\\
              (1) Faculty of Natural Sciences and Mathematics, University of Maribor, Koro\v{s}ka 160, SI-2000 Maribor,
   Slovenia; \\(2) Institute of Mathematics, Physics and Mechanics, Jadranska 19, SI-1000 Ljubljana, 
   Slovenia; \\(3) Andrej Maru\v si\v c Institute, University of Primorska, Muzejski trg 2, SI-6000 Koper,
   Slovenia\\
             {iztok.banic@um.si}           %  \\
     
				\-
				
		\noindent G.  Erceg\\
             Faculty of Science, University of Split, Rudera Bo\v skovi\' ca 33, Split,  Croatia\\
%             {i}     
{{goran.erceg@pmfst.hr}       }    %  \\

                 	\-
                 				
		\noindent A.  Illanes\\
             Instituto de Matemáticas, Universidad Nacional Autónoma de México, Circuito Exterior, Ciudad Universitaria, C.P. 04510, Ciudad de México, México\\
%             {i}     
{{illanes@ciencias.unam.mx}       }    %  \\

           \-
                 				
		\noindent I.  Jeli\' c\\
             Faculty of Science, University of Split, Rudera Bo\v skovi\' ca 33, Split,  Croatia\\
%             {i}     
{{ivajel@pmfst.hr}       }    %  \\

                 	\-
					
  \noindent J.  Kennedy\\
             Department of Mathematics,  Lamar University, 200 Lucas Building, P.O. Box 10047, Beaumont, Texas 77710 USA\\
%             {}     
{{kennedy9905@gmail.com}       }

	\-
				
		\noindent V.  Nall\\
             Department of Mathematics,  University of Richmond, Richmond, Virginia 23173 USA\\
%             {i}     
{{vnall@richmond.edu}       }      %  \\


                 %  \\

%``text''
%%%%%%%%%%%%%%%%%%%%%%%%%%%%%%%%%%%%%%%%%%%%%%%%%%%%%%%%%%%%%%%%%%%%%%%%%%%%%%%%%
%%% I N T R O D U C T I O N S

\end{document}